**Miquel point and isogonal conjugation**
Valcho Milchev

**Abstract.** We study the possible positions of the Miquel point in the plane of a given triangle, which Miquel triangles are similar to the given one. We found out that these positions are eleven. We also study the possible positions of the Miquel point in the plane of a given triangle, where among the families of Miquel triangles there are triangles, which are similar to the given triangle. We study which of them are isogonal conjugated.

## 1. Introduction

First we give some facts from the theory of Miquel points.

Let $ABC$ be a triangle with angles $\sphericalangle A = \sphericalangle BAC$, $\sphericalangle B = \sphericalangle CBA$ and $\sphericalangle C = \sphericalangle ACB$ respectively.

**Theorem 1.** (Miquel Theorem) $[1, p.79]$ **If three points** $X$, $Y$ and $Z$ **are chosen on each side of triangle** $ABC$ **or on their extensions, then the three circles through any vertex of the triangle** $\triangle ABC$ **and two points on the sides of the triangle intersect at one point** (Figure 1a, 1b, 1c).

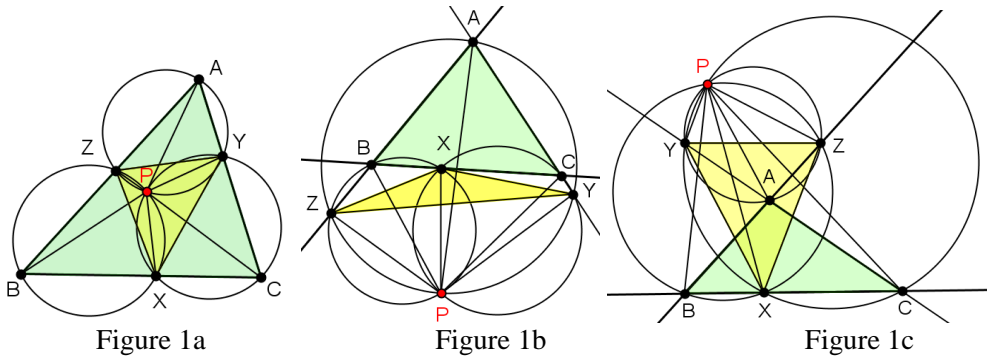

Figure 1a  Figure 1b  Figure 1c

**Definition 1**. The point $P$ from Theorem 1 is called ***Miquel point,*** and the three circles are called Miquel circles.

**Definition 2**. The triangle $XYZ$ is called ***Miquel triangle*** of the point $P$ relative to triangle $\triangle ABC$, and the three vertexes of triangle $\triangle XYZ$ are called Miquel triad $(X,Y,Z)$.

**Remark 1.** In the special case, when point $P$ is on the circum-circle, the triple of points $(X,Y,Z)$ is collinear (***generalized Simson Theorem)***.

**Remark 2.** We note that from the Miquel Theorem it follows that any point in the plane of a given triangle can be considered as Miquel point. We also denote such a set of triple of points with $(X,Y,Z)$ - that set is called family $(X,Y,Z)$ of Miquel triangles.

We denote the interior angles of triangle $\triangle XYZ$ with
$\sphericalangle X = \sphericalangle YXZ$, $\sphericalangle Y = \sphericalangle ZYX$ and $\sphericalangle Z = \sphericalangle XZY$, respectively.

**Theorem 2.** (***Miquel equations***) $[1, p.84]$ **For the Miquel point are in force these equations with directed angles:**

$$\angle A + \angle X = \angle BPC, \quad \angle B + \angle Y = \angle CPA, \quad \angle A + \angle Z = \angle APB.$$

**Corollary 1.** All Miquel triangles of any fixed point $P$ relative to triangle $\triangle ABC$ are similar to themselves. One of these triangles is the pedal triangle of point $P$.

**2. Points which Miquel triangles are similar to the given one**

Given the basic triangle $ABC$. In this section we study the positions of point $P$, which Miquel triangle $XYZ$ is similar to the given one. We make a thorough study of non-trivial cases – for scalene triangle and non-right triangle $ABC$, and point $P$ is not on any of the line $AB$, $BC$ and $CA$.

In the process of studying the positions of Miquel point, which Miquel triangles are similar to the given one, we use the following proposition.

**Theorem 3.** $[2, p.43]$ Two points, which are inverse with respect to the circum-circle of a given triangle, have similar pedal triangles.

This means that finding all Miquel points, which triands are triangles similar to the given, will be comprehensively if all Miquel points in the circum-circle of triangle $ABC$ are found.

Therefore we will consider the positions of the Miquel point in the circum-circle. We will use two important lemmas.

**Lemma 1.**

(i) If the Miquel point $P$ is interior for triangle $ABC$, then point $P$ is interior for the Miquel triangle $XYZ$ too.

(ii) If the Miquel point $P$ is exterior for triangle $ABC$, then point $P$ is exterior for the Miquel triangle $XYZ$ too.

**Proof.** If the Miquel point $P$ is interior for triangle $\triangle ABC$, then $\angle XPY + \angle YPZ + \angle ZPX = 360°$. This is because the three angles in the left side of this equation supplements the interior angles $\angle A$, $\angle B$ и $\angle C$, respectively (Figure 1a). Consequently point $P$ is also interior for the Miquel triangle $\triangle XYZ$.

If the Miquel point $P$ is exterior for triangle $\triangle ABC$ (Figure 1b, 1c), then the inequality $\angle XPY + \angle ZPY + \angle ZPX < 360°$ is in force. Here two of the angles are equal to two of the interior angles, and the third one supplements the rest interior angle, for example

$$\angle XPY + \angle ZPY + \angle ZPX = \angle C + (180° - \angle A) + \angle B =$$
$$= \angle C + (180° - \angle A) + \angle B = \angle A + \angle B + \angle C + 180° - 2\angle A = 360° - 2\angle A < 360°.$$

Consequently in this case point $P$ is also exterior for the Miquel triangle $\triangle XYZ$.

For the sake of brevity and convenience we will give some symbols:

Given triangle $ABC$ and point $P$. Then denote with directed angles (Figure 2a, 2b):
$\alpha_1 = \angle PAC$, $\alpha_2 = \angle BAP$, $\beta_1 = \angle PBA$, $\beta_2 = \angle CBP$, $\gamma_1 = \angle PCB$, $\gamma_2 = \angle ACP$.

**Lemma 2.** Given point $P$ interior for the circum-circle of triangle $ABC$ and $XYZ$ is Miquel triangle of point $P$ relative to triangle $\triangle ABC$. Then the following equations for the angles of triangle $\triangle XYZ$ are in force:

$\sphericalangle X = \beta_1 + \gamma_2$, $\sphericalangle Y = \gamma_1 + \alpha_2$, $\sphericalangle Z = \alpha_1 + \beta_2$.

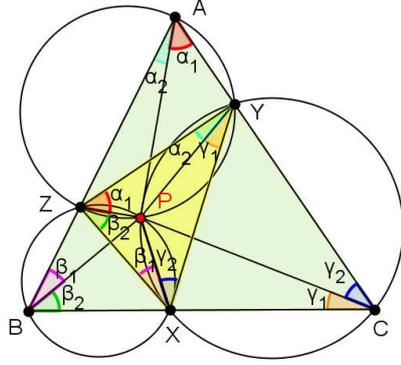
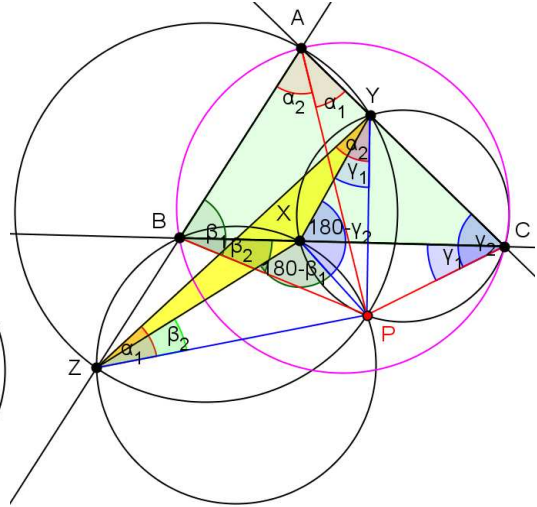

Figure 2a          Figure 2b

**Proof**. Equations with directed angles we prove using inscribed angles with the same adjacent arc. In Figure 2a is the case where the point $P$ is in the interior of triangle $\triangle ABC$. In Figure 2b is the case where point $P$ is in the exterior of triangle $\triangle ABC$, but in the circum-circle of triangle $ABC$ - now $\beta_2$ and $\gamma_1$ have negative digits.

**Theorem 4. In the plane of a given triangle** $ABC$ there are 11 points, which Miquel triangles relative to triangle $\triangle ABC$ are similar to the given one. Six of them are interior for the circum-circle of the given triangle and the rest five points are exterior for the circum-circle. The interior points are First Brocard point, Second Brocard point, the center of the circum-circle and three points on the symmedians $S_A$, $S_B$ and $S_C$ of triangle $\triangle ABC$. The exterior points are the inverse points of First Brocard point, Second Brocard point and the three points on the symmedians - $S_A$, $S_B$ and $S_C$ ($S_A$ is the point of intersection of the ray of the symmedian from vertex $A$ and the arc $\overarc{BOC}$, likewise $S_B$ and $S_C$) with respect to the circum-circle.

**Proof.** Let list the possible similarities of $\triangle ABC$ and $\triangle XYZ$ - there are 6 combinations: $\triangle ABC \sim \triangle XYZ$; $\triangle ABC \sim \triangle ZXY$; $\triangle ABC \sim \triangle YZX$; $\triangle ABC \sim \triangle XZY$; $\triangle ABC \sim \triangle ZYX$; $\triangle ABC \sim \triangle YXZ$. Similarities are listed in a common order – the corresponding angles are equal, for example $\triangle ABC \sim \triangle ZXY$ means that $\sphericalangle Z = \sphericalangle A$, $\sphericalangle X = \sphericalangle B$, $\sphericalangle Y = \sphericalangle C$.

We consider the similarities separately applying Lemma 2.

**1.** Let $\triangle ABC \sim \triangle XYZ$. That means $\sphericalangle X = \sphericalangle A$, $\sphericalangle Y = \sphericalangle B$, $\sphericalangle Z = \sphericalangle C$. From here and the Miquel equation (Theorem 2) it follows that $\sphericalangle BPC = \sphericalangle A + \sphericalangle X = 2\sphericalangle A$. Similarly we prove $\sphericalangle CPA = \sphericalangle B + \sphericalangle Y = 2\sphericalangle B$, $\sphericalangle APB = \sphericalangle C + \sphericalangle Z = 2\sphericalangle C$. These equations determine that in this case point $P$ coincides with point $O$ - the center of the circum-circle of triangle $\triangle ABC$.

**2**. Let $\triangle ABC \sim \triangle ZXY$. In this case $\sphericalangle X = \sphericalangle B$, $\sphericalangle Y = \sphericalangle C$, $\sphericalangle Z = \sphericalangle A$. Consequently $\beta_1 + \gamma_2 = \beta_1 + \beta_2$, $\gamma_1 + \alpha_2 = \gamma_1 + \gamma_2$, $\alpha_1 + \beta_2 = \alpha_1 + \alpha_2$, whence $\alpha_2 = \beta_2 = \gamma_2$ and now point $P$ coincides with the First Brocard point of triangle $\triangle ABC$.

**3**. Let $\triangle ABC \sim \triangle YZX$. Now $\sphericalangle X = \sphericalangle C$, $\sphericalangle Y = \sphericalangle A$, $\sphericalangle Z = \sphericalangle B$.

Consequently $\beta_1 + \gamma_2 = \gamma_1 + \gamma_2$, $\gamma_1 + \alpha_2 = \alpha_1 + \alpha_2$, $\alpha_1 + \beta_2 = \beta_1 + \beta_2$, whence $\alpha_1 = \beta_1 = \gamma_1$ - here point $P$ coincides with the Second Brocard point of triangle $\triangle ABC$.

In the next three cases we use some facts for the symmedians in a triangle.

**Definition 3**. Symmedian of a triangle is a line that is symmetrical to the median with respect to the angle bisector of the interior angle of the triangle.

**Lemma 3.** $[1, p.58]$ Let point $D$ is on the side $BC$ of triangle $ABC$. The line $AD$ contains the symmedian from vertex $A$ if and only if point $D$ divides the side $BC$ into parts that are proportional to the squares of the adjacent sides, i.e. $\dfrac{BD}{CD} = \left(\dfrac{AB}{AC}\right)^2$.

**4. Now we consider the case** $\triangle ABC \sim \triangle XZY$. We will prove that point $P$ is a

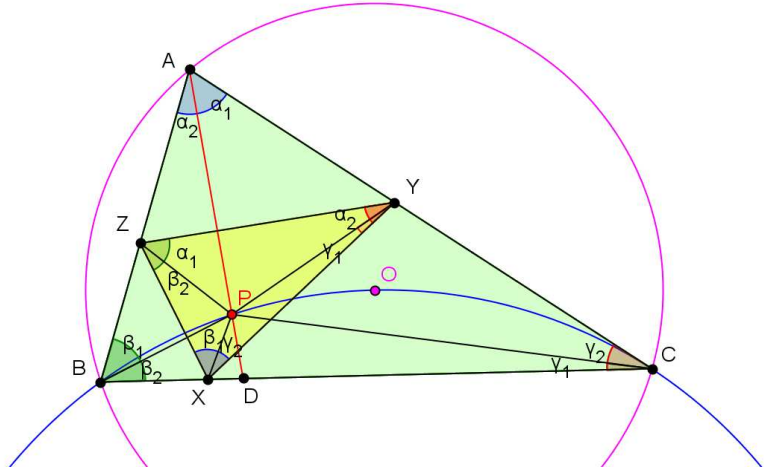

Figure 3

point of intersection of the ray of the symmedian from vertex $A$ and the arc $\widehat{BOC}$ from the circle through points $B$, $C$ and the center $O$ of the circum-circle of triangle $\triangle ABC$. In Figure 3 it is shown the case $\angle A < 90°$, and in Figure 4 - $\angle A > 90°$.

First we consider the case $\angle A < 90°$. Since $\triangle ABC \sim \triangle XZY$, then $\angle X = \angle A$, $\angle Y = \angle C$, $\angle Z = \angle B$. According to Lemma 2 we have $\beta_1 + \gamma_2 = \alpha_1 + \alpha_2$, $\gamma_1 + \alpha_2 = \gamma_1 + \gamma_2$, $\alpha_1 + \beta_2 = \beta_1 + \beta_2$, whence $\gamma_2 = \alpha_2$ and $\beta_1 = \alpha_1$. Applying the Miquel equation (Theorem 2) we obtain
$\angle BPC = \angle A + \angle X = 2\angle A$,
$\angle CPA = \angle B + \angle Y = \angle B + \angle C = 180° - \angle A$,
$\angle APB = \angle C + \angle Z = \angle C + \angle B = 180° - \angle A$.

Consequently, if exists, point $P$ is on the arc from which the side $BC$ is visible at angle $2\angle A$, then point $P$ is also on the arc from which the side $CA$ is visible ate angle $180° - \angle A$, and point $P$ is on the arc from which the side $AB$ is visible at angle $180° - \angle A$. Such a point, if there is any at all, is unique in the interior of the circum-circle of triangle $\triangle ABC$. Let check it.

Let $AP$ intersects $BC$ at point $D$. Then $PD$ is an angle bisector in triangle $\triangle BPC$ and consequently $\dfrac{BD}{DC} = \dfrac{BP}{PC}$. But $\triangle BPA \sim \triangle APC$, whence $\dfrac{AB}{AC} = \dfrac{BP}{AC} = \dfrac{AC}{CP}$ and $\left(\dfrac{AB}{AC}\right)^2 = \dfrac{BP}{AC} \cdot \dfrac{AC}{CP} = \dfrac{BP}{CP} = \dfrac{BD}{DC}$. The last equation shows that point $P$ exists and the symmedian from vertex $A$ of triangle $\triangle ABC$ contains point $P$.

In the case $\sphericalangle A > 90°$ (Figure 4), we obtain the same result. Now point $P$ is exterior for triangle $\triangle ABC$, but it is also on the arc $\overset{\frown}{BOC}$ from the circle through points $B$, $C$ and the center $O$ of the circum-circle of triangle $\triangle ABC$, and the line $PA$ contains the symmedian from vertex $A$ and it is an angle bisector in triangle $\triangle BPC$.

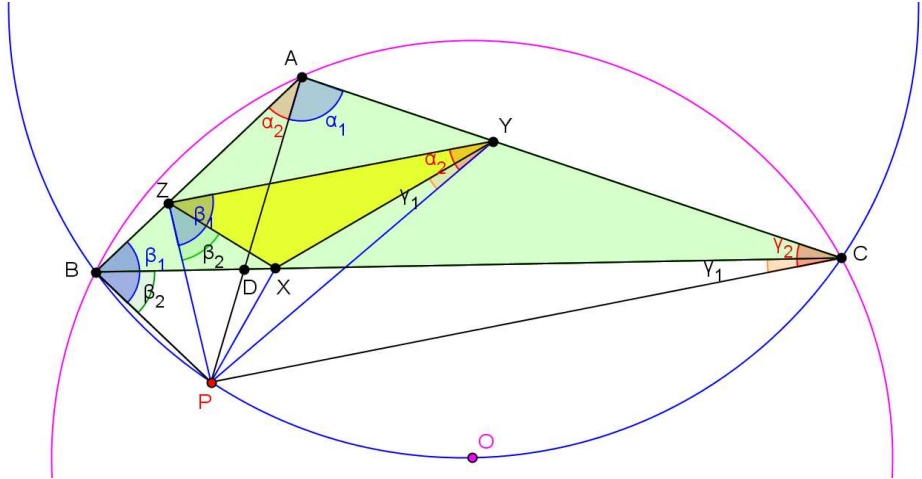

Figure 4

Consequently the demand point is the intersection of the symmedian from vertex $A$ and the arc through points $B$ and $C$, which contains the center of the circum-circle. Let denote this point with $S_A$.

Similarly in the cases $\triangle ABC \sim \triangle ZYX$ and $\triangle ABC \sim \triangle YXZ$ we find out points $S_B$ and $S_C$, which are on the symmedians from vertexes $B$ and $C$, respectively.

Consequently there are exactly 6 points in the interior of a given triangle $\triangle ABC$, which Miquel triangles are similar to the given one – these points are first Brocard point, second Brocard point, the center of the circum-circle and the points $S_A$, $S_B$ and $S_C$ on the symmedians from vertexes $A$, $B$ and $C$, respectively.

Taking account of that for each point $P$ its pedal triangle is one of the family of triangles $\triangle XYZ$, and the pedal triangles of the inverse points with respect to the circum-circle are similar, we some to the conclusion that in the plane of a given triangle there are 5 points, which Miquel triangles are similar to the given one – these points are the inverse points with respect to the circum-circle of the first Brocard point, the second Brocard point and the points $S_A$, $S_B$ and $S_C$, since the inverse point of the center of the circum-circle is an infinite point.

### 3. Isogonal conjugated points of Miquel point

We will consider concrete positions of the Miquel point starting with the positions for which Miquel triangles are similar to the main triangle.

**Theorem 5.** If Miquel point $P$ is center of the circum-circle of the main triangle $\triangle ABC$, then point $P$ is orthocenter of the Miquel triangle $\triangle XYZ$.

**Proof**. If point $P$ is center of the circum-circle of triangle $\triangle ABC$, then $\alpha_2 = \beta_1$, $\beta_2 = \gamma_1$, $\gamma_2 = \alpha_1$ (Figure 5a, 5b). We obtain the equation $\gamma_2 + \gamma_1 + \beta_1 = \alpha_1 + \beta_2 + \alpha_2 = 90°$. That means that $YP \perp ZX$. Similarly $ZP \perp XY$ and $XP \perp YZ$. Consequently point $P$ is orthocenter of triangle $\triangle XYZ$.

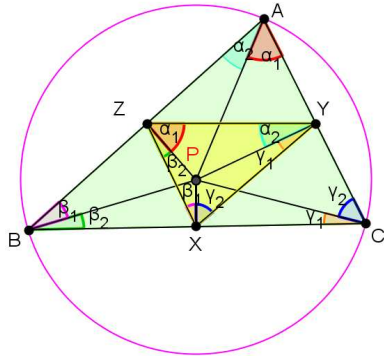
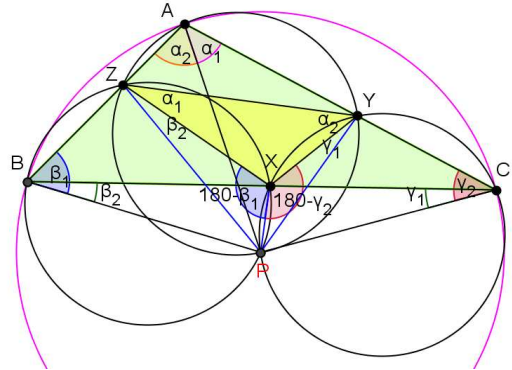

Figure 5a                              Figure 5b

**Theorem 6.** If Miquel point $P$ is orthocenter of triangle $\triangle ABC$, then point $P$ is the center of the inscribed circle or the excircle of the Miquel triangle $\triangle XYZ$ (Figure 6a, 6b).

**Proof**. We will consider two cases:

(1) Let point $P$ is orthocenter of the acute-angled triangle $\triangle ABC$ (Figure 6a). Then $\gamma_1 = \alpha_2$, $\alpha_1 = \beta_2$, hence in triangle $\triangle XYZ$ point $P$ is the intersection of the interior angle bisectors, i.e. it is the center of the inscribed circle.

(2) Let $\triangle ABC$ is an obtuse-angled triangle and let $\sphericalangle A > 90°$ (Figure 6b). In this case $\beta_1 = \gamma_2$, $\gamma_1 = 180° - \alpha_2$, $\alpha_1 = 180° - \beta_2$, hence $PY$ and $PZ$ are exterior angle bisectors, and $PX$ is an interior angle bisector. Consequently point $P$ is center of the excircle of triangle $\triangle XYZ$.

**Theorem 7.** If miquel point $P$ is center of the incircle or excircle of triangle $\triangle ABC$, then point $P$ is center of the circum-circle of triangle $\triangle XYZ$.

**Proof**. We consider two cases (Figure 7a, 7b):

(1) Let $P$ is center of the incircle of triangle $\triangle ABC$. Then $\alpha_1 = \alpha_2$, $\beta_1 = \beta_2$, $\gamma_1 = \gamma_2$, hence point $P$ is on the bisectors of the sides of triangle $\triangle XYZ$ and therefore it is center of the circum-circle of triangle $\triangle XYZ$.

(2) Let $P$ is center of the excircle to side $BC$ of triangle $\triangle ABC$. And in this case , $PB$ and $PC$ are bisectors of the sides of triangle $\triangle XYZ$.

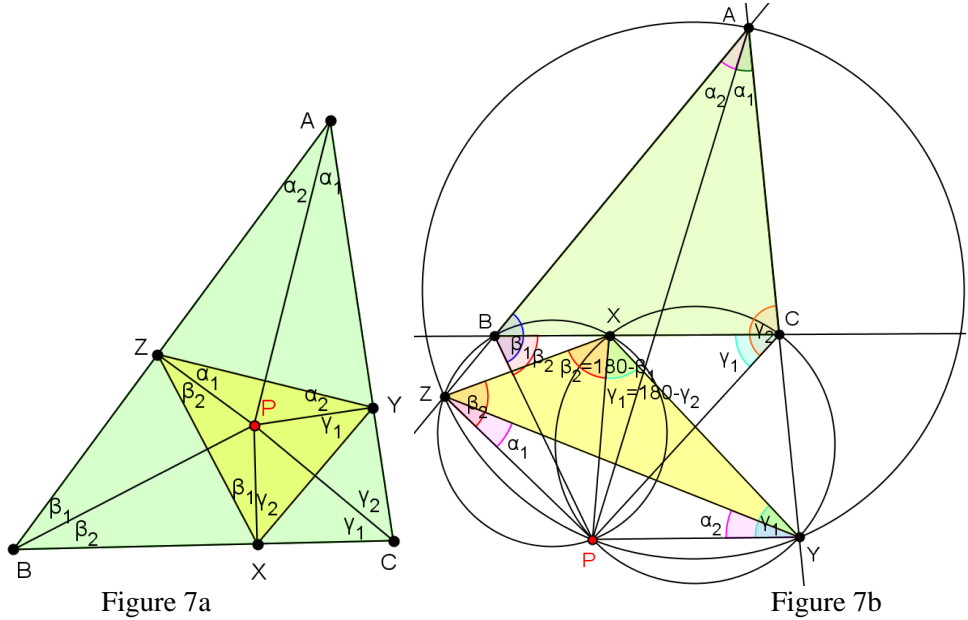

Figure 7a        Figure 7b

**Theorem 8.**

(i) If Miquel point $P$ is first Brocard point of triangle $\triangle ABC$, then it is the also first Brocard point of the Miquel triangle $XYZ$.

(ii) If Miquel point $P$ is second Brocard point of triangle $\triangle ABC$, then it is also second Brocard point of the Miquel triangle $XYZ$.

**Proof.** We will use the symbols in figure 2a.

(i) If point $P$ is first Brocard point of triangle $\triangle ABC$, then the equation $\alpha_2 = \beta_2 = \gamma_2$ is fulfilled, whence from Lemma 2 it follows that point $P$ is first Brocard point of triangle $\triangle XYZ$.

(ii) If point $P$ is second Brocard point of triangle $\triangle ABC$, then the equation $\alpha_1 = \beta_1 = \gamma_1$ is fulfilled and from Lemma 2 it follows that point $P$ is second Brocard point of triangle $\triangle XYZ$.

Now we will consider the Miquel triangles of points $S_A$, $S_B$ and $S_C$.

**Theorem 9.** Let Miquel point $P$ coincides with point $S_A$ of triangle $\triangle ABC$. Then in the Miquel triangle $XYZ$ point $P$ is point $M_X$ on the median $XE$ from vertex $X$, for which $M_X E = EF$, if $\sphericalangle A < 90°$, and $XE = EF$, if $\sphericalangle A > 90°$, where point $E$ is the middle of $YZ$, and point $F$ is the intersection point of $XM_x$ and the circum-circle of triangle $\triangle XYZ$.

**Proof.** Let $\triangle ABC \sim \triangle XZY$, i.e. $\sphericalangle X = \sphericalangle A$, $\sphericalangle Y = \sphericalangle C$, $\sphericalangle Z = \sphericalangle B$. We will use the received symbols of the angles of the triangles and we will consider two cases – if $\sphericalangle A < 90°$ and if $\sphericalangle A > 90°$.

(1) If $\sphericalangle A < 90°$ (Figure 8). Then $\sphericalangle XPY = 180° - (\gamma_1 + \gamma_2) = 180° - \sphericalangle C$, $\sphericalangle ZPX = 180° - (\beta_1 + \beta_2) = 180° - \sphericalangle B$, hence $\sphericalangle YPF = \sphericalangle C = \sphericalangle Y$, $\sphericalangle FPZ = \sphericalangle B = \sphericalangle Z$. From the equations of inscribed angles It follows that $\sphericalangle PFY = \sphericalangle Z = \sphericalangle B = \sphericalangle FPZ$ and $\sphericalangle YPF = \sphericalangle Y = \sphericalangle C = \sphericalangle ZFP$. That means that the quadrangle $PYFZ$ is a parallelogram and $XE$ is a median of triangle $\triangle XYZ$.

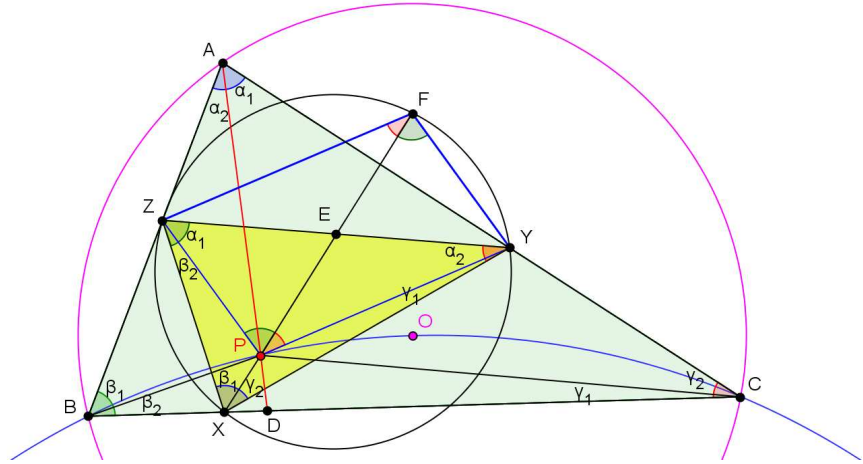

Figure 8. Here $P \equiv S_A \equiv M_X$.

(2) If $\sphericalangle A > 90°$ (Figure 9). We extend $XP$ till it intersects the Miquel circle through vertex $A$ in point $F$. We have $\alpha_1 = \beta_1$, $\alpha_2 = \gamma_2$, because $P \equiv S_A$. Then $XYFZ$ is a parallelogram and $PF$ contains the median from vertex $X$ in triangle $\triangle XYZ$. Again $XE = EF$ and $ZE = EY$. Here $\sphericalangle FXZ = \sphericalangle PFY = \sphericalangle PAC = \alpha_1$, because point $P$ is Miquel point of triangle $\triangle ABC$ and point on the symmedian $S_A$ at the same time. On the other hand $\sphericalangle PFY = \sphericalangle PAC = \alpha_1$, as inscribed angles with a same

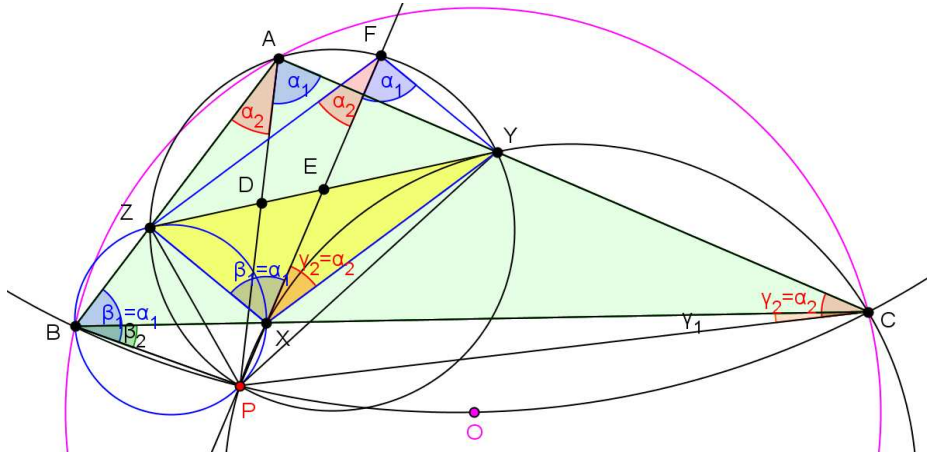

Figure 9. Here $P \equiv S_A \equiv M_X$.

arc.

**Theorem 10.** Let Miquel point $P$ coincides with point $M_A$ of triangle $\triangle ABC$. Then the Miquel triangle $\triangle XYZ$ is isosceles one, point $P$ is on the arc $\widehat{YLZ}$, which contains the intersection point $L$ of the interior angle bisectors of triangle $\triangle XYZ$ /the center of the incircle/, and point $P$ is an interior point of triangle $\triangle XYZ$, if $\triangle ABC$ is acute-angled and it is exterior point of triangle $\triangle XYZ$, if it is obtuse-angled.

**Proof**. We also consider two cases.

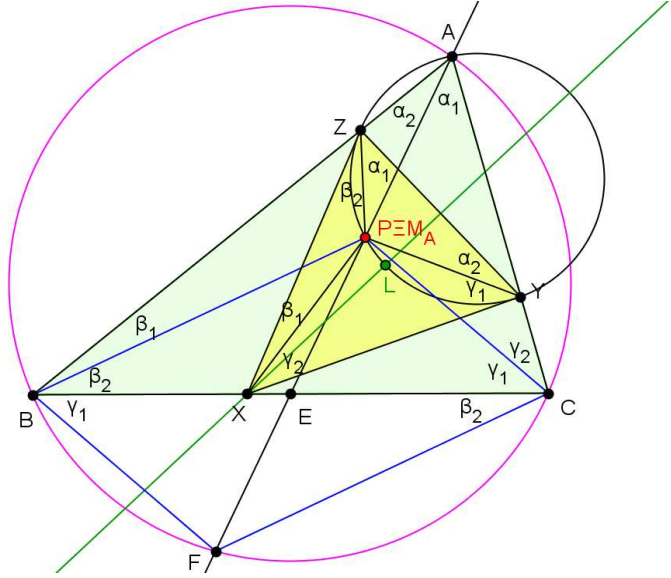

Figure 10. Here $P \equiv M_A$.

(1) If triangle $\triangle ABC$ is acute-angled (Figure 10). From the separation of the angles of Miquel point relative to $\triangle ABC$ it is clear that in triangle $\triangle XYZ$ we have $\angle ZYX = \angle XZY = \angle A$ and $\angle YPZ = 180° - \angle A$, which means that point $P$ „sees" the side $YZ$ at $180° - \angle A$. We come to the conclusion that point $P$ is on the arc $\widehat{ZLY}$, where $L$ is an intersection point of the interior angle bisectors of triangle $\triangle XYZ$ (the center of the incircle). The equations $\angle ZYP = \angle XZP = \alpha_2$ and $\angle PZY = \angle PYX = \alpha_1$ show us that $XY$ and $XZ$ are tangents of the circle in points $Y$, $Z$, $L$ and $P$. That circle is one of the Miquel circle of triangle $\triangle ABC$.

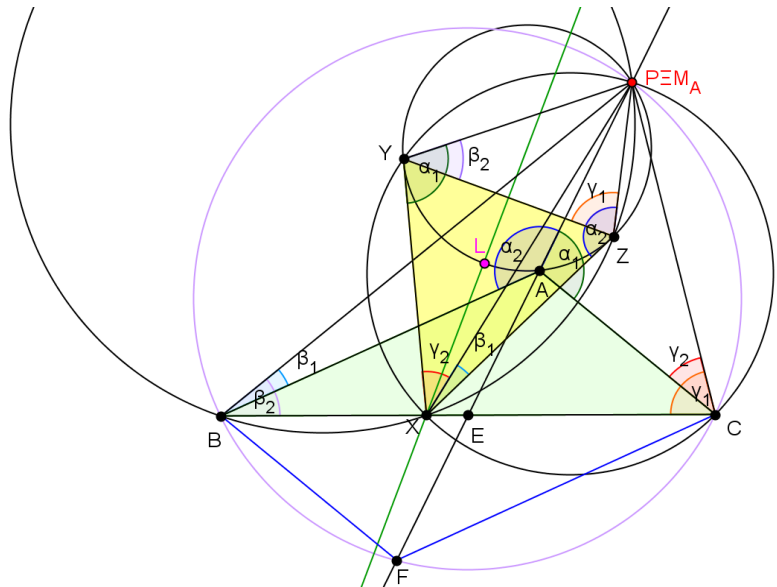

Figure 11. Here $P \equiv M_A$.

(2) If the triangle is obtuse-angled one– let $\angle A > 90°$ (Figure 11). In this case point $P$ is exterior for triangle $\triangle ABC$, it is on the arc from the circle through points $L$, $Y$ and $Z$, which doesn't contain point $L$ and from it $YZ$ is visible at $180° - \angle A$. That circle pass throughout the intersection point of the interior angle bisectors of triangle $\triangle ABC$ and it is one of the Miquel circles of triangle $\triangle ABC$ wit Miquel point $P$. The sides $XY$ and $XZ$ are tangents of the circle in points $Y$, $Z$, $L$ and $P$. Indeed in this case there are the following equations between angles: $\angle XPZ = \angle B$, $\angle YPX = \angle C$, $\angle YPZ = \angle B + \angle C = 180° - \angle A$ and $\angle ZAL = \angle ZYL = \angle A$.

**Theorem 11.** Let $\triangle ABC$ is an isosceles triangle ($AB = AC$) and Miquel point $P$ is on the arc $\overset{\frown}{CLB}$, where $L$ is the intersection point of the interior angle bisectors. Then point $P$ coincides with point $S_X$ in the Miquel triangle $\triangle XYZ$. ($S_X$ is the intersection point of the symmedian from vertex $X$ and the arc $\overset{\frown}{YOZ}$ in $\triangle XYZ$, whereo $O$ is center of the circum-circle).

**Proof**.

(1) If $\angle A < 90°$. In the scheme in Figure 12a it is shown that the angles of triangle $\triangle XYZ$ are connected with the equations $\angle YPZ = \angle B + \angle C = 2\angle X$, $\angle XPY = \angle ZPX = 180° - \angle X$ and point $P$ coincides with $S_X$ in $\triangle XYZ$.

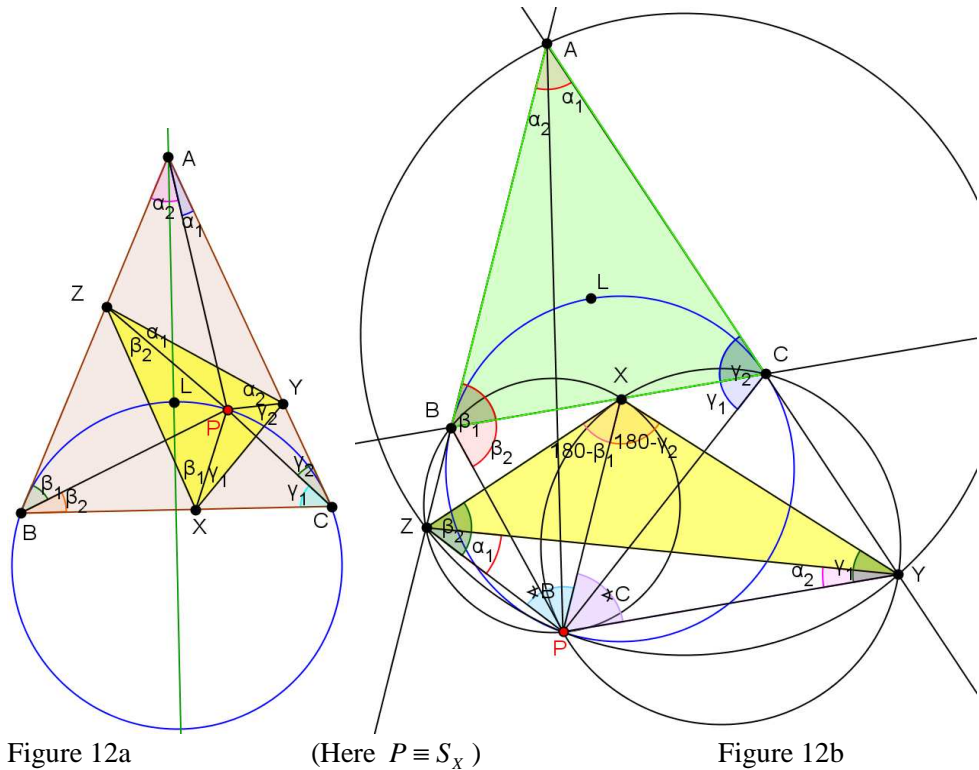

Figure 12a        (Here $P \equiv S_X$)        Figure 12b

(2) If $\angle A > 90°$. In this case point $P$ coincides with $S_X$ in $\triangle XYZ$, which is shown in Figure 12b.

Between the considered positions of Miquel point there is a remarkable connection – there is an isogonal conjugation of some couple of points.

**Definition 4**. Points $P$ and $Q$ are isogonal conjugated relative to triangle $ABC$, if $\angle CBQ = \angle PBA$, $\angle BAP = \angle QAC$, $\angle ACQ = \angle PCB$ (Figure 13abc).

Between the considered points in this paper there are points, for which there is known isogonal conjugation [3]:

-Well-known points that are isogonal conjugated are the first Brocard point and the second Brocard point.
-The center of the circum-circle and the orthocenter are isogonal conjugated.
-The center of the incircle and the centers of the excircles are isogonal conjugated along.

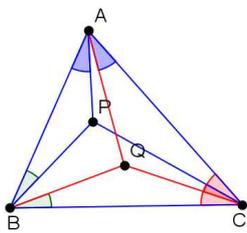
Figure 13a

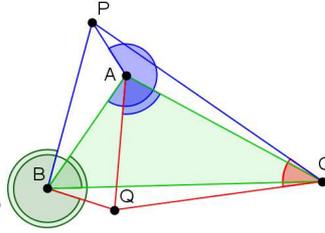
Figure 13b

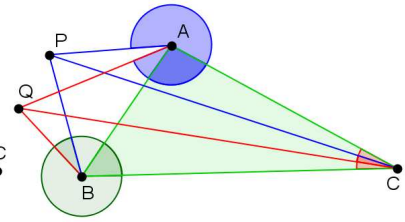
Figure 13c

**Theorem 12.** Let $AE$ and $AD$ are median and symmedian of triangle $ABC$, respectively, where points $D$ and $E$ are on the side $BC$. Let $S_A$ is the unique point, for which the condition $\sphericalangle BS_AC = 2\sphericalangle A$, $\sphericalangle CS_AA = \sphericalangle AS_AB = 180° - \sphericalangle A$ is fulfilled. Let points $M_A$ and $F$ are on the median $AE$ of triangle $ABC$ such that

(i) if $\sphericalangle A < 90°$ the median $AE$ intersects the circum-circle of triangle $ABC$ for second time in point $F$ and $EM_A = EF$ is fulfilled (Figure 14a).

(ii) if $\sphericalangle A > 90°$ the quadrangle $ABFC$ is a parallelogram, and the median $AE$ intersects the circum-circle of triangle $FBC$ for second time in point $M_A$ (Figure 14b).

Then points $S_A$ and $M_A$ are isogonal conjugated relative to $\triangle ABC$.

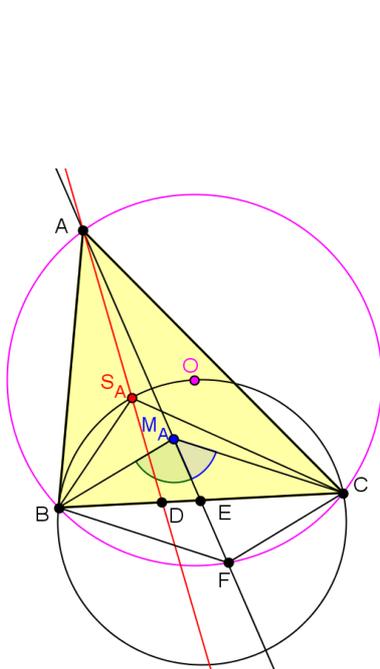
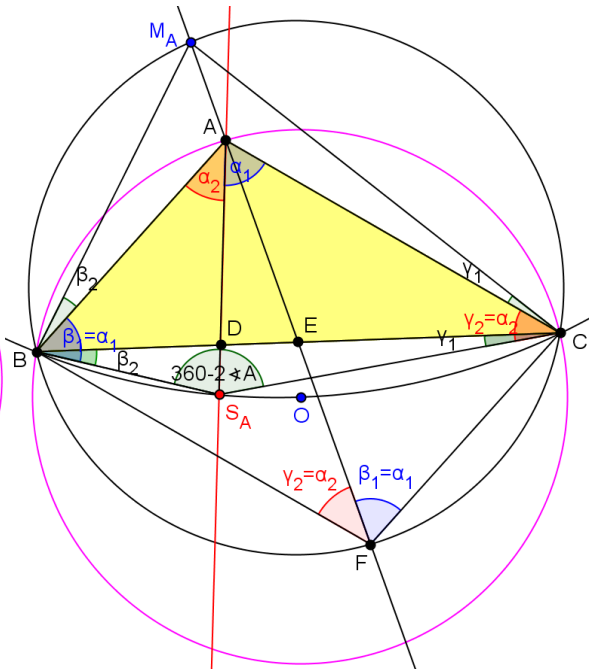

Figure 14a  Figure 14b

**Proof.** We will consider two cases – if $\sphericalangle A < 90°$ and if $\sphericalangle A > 90°$.

(i) Let $\sphericalangle A < 90°$. We have [×], because $AM_A$ and $AS_A$ are median and symmedian of triangle $ABC$. Hence $\sphericalangle BAM_A = \sphericalangle S_AAC$. On the other hand $S_A$ is a special point on the symmedian from vertex $A$, for which $\sphericalangle S_A AC = \sphericalangle S_A BA$. Furthermore from the properties of the angles in a parallelogram and the inscribed angles we obtain the equations $\sphericalangle BAM_A = \sphericalangle FCB = \sphericalangle M_A BC$. We conclude $\sphericalangle CBM_A = \sphericalangle S_A BA$. Similarly we obtain the equation $\sphericalangle M_A CB = \sphericalangle ACS_A$. Consequently points $S_A$ and $M_A$ are isogonal conjugated relative to triangle $ABC$.

(ii) Let $\sphericalangle A > 90°$. Similarly we come to conclusion that points $S_A$ and $M_A$ are isogonal conjugated relative to triangle $\triangle ABC$.

**Theorem 13.** Let in an isosceles triangle $ABC$, $AB = AC$, point $Q$ is on the circle through points $B$, $C$ and $L$, where $L$ is the intersection point of the interior angle bisectors of triangle $\triangle ABC$ (center of the incircle), and point $T$ is symmetrical of point $Q$ relative to the interior angle bisector from vertex $A$. Then $Q$ and $T$ are isogonal conjugated points relative to triangle $ABC$ (Figure 15a, 15b).

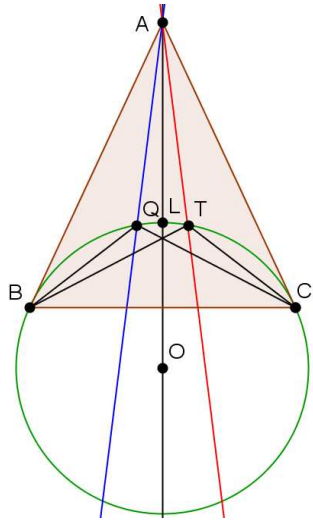 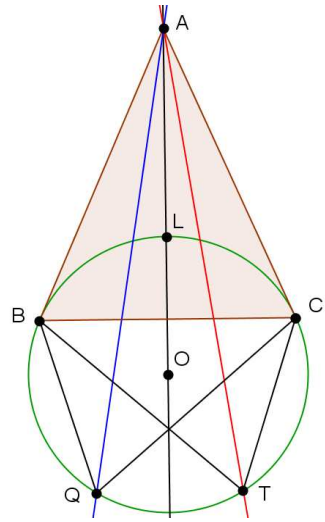

Figure 15a        Figure 15b

**Proof.** In this case $AB$ and $AC$ are tangents of the circle through points $B$, $C$ and $L$. Consequently $\sphericalangle ACT = \sphericalangle QBA$, $\sphericalangle CBT = \sphericalangle QCB$ and together with the symmetry of $AQ$ and $AT$ we obtain that points $Q$ and $T$ are isogonal conjugated relative to triangle $ABC$.

### 4. Chains of Miquel triads

Further we consider families of Miquel points relative to a given main triangle with fixed Miquel point, which are obtained as follows:

In the plane are given the main triangle $\triangle A_0 B_0 C_0$ and a fixed point $P$. Then triangles are constructed, as on each of the sides of $\triangle A_k B_k C_k$ or on their continuations is taken one point - $A_{k+1}$ on the side $B_k C_k$, $B_{k+1}$ on the side $C_k A_k$ and $C_{k+1}$ on the side $A_k B_k$ in such a way that $\triangle A_{k+1} B_{k+1} C_{k+1}$ is a Miquel triangle of $\triangle A_k B_k C_k$ and point $P$ with interior angles $\sphericalangle A_k = \sphericalangle B_k A_k C_k$, $\sphericalangle B_k = \sphericalangle C_k B_k A_k$ and $\sphericalangle C_k = \sphericalangle A_k C_k B_k$, where $k = 0,1,2,3,...$. So $\triangle A_{k+1} B_{k+1} C_{k+1}$ is a Miquel triangle of $\triangle A_k B_k C_k$ and point $P$. In each of these families of Miquel triangles one is the pedal triangle of point $P$ relative to the corresponding main triangle.

All triangle from a given family are countless and they are similar to each other. In Figure 16 there are triangles from the families $(A_1, B_1, C_1)$, $(A_2, B_2, C_2)$, $(A_3, B_3, C_3)$.

In the queue of Miquel triangle with a fixed Miquel point there are periodically repeating similar triangles.

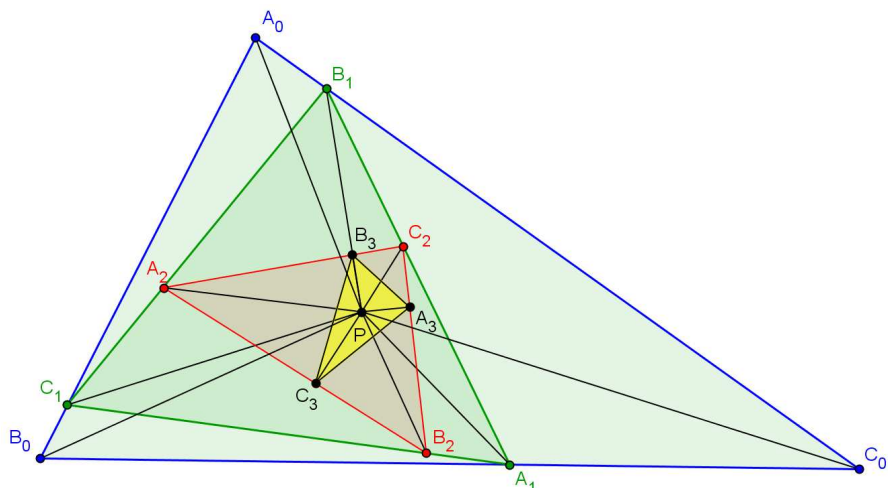

Figure 16

**Theorem 14.** [4] For the queue of Miquel triangles $\{\triangle A_k B_k C_k\}$, where $k = 0, 1, 2, 3, \ldots$, is in force $\triangle A_i B_i C_i \sim \triangle A_j B_j C_j$, if $i \equiv j \pmod{3}$.

The obtained results in this paper give us chance to list all similarities in the queues of families of Miquel triads.

**Theorem 15.**

- if Miquel $P$ is also first Brocard point of triangle $\triangle A_0 B_0 C_0$, then it is first Brocard point of Miquel triangles $A_k B_k C_k$ and $\triangle A_0 B_0 C_0 \sim \triangle A_k B_k C_k$ for each $k = 1, 2, 3, \ldots$.

- if Miquel point $P$ is also second Brocard point of triangle $\triangle A_0 B_0 C_0$, then it is second Brocard point of Miquel triangles $A_k B_k C_k$ and $\triangle A_0 B_0 C_0 \sim \triangle A_k B_k C_k$ for each $k = 1, 2, 3, \ldots$.

- if Miquel point $P$ is center of the circum-circle or one of the points $S_{A0}$, $S_{B0}$, $S_{C0}$ of triangle $\triangle A_0 B_0 C_0$, then $\triangle A_0 B_0 C_0 \sim \triangle A_k B_k C_k$ for each $k \equiv 0 \pmod{3}$ and $k \equiv 1 \pmod{3}$.

- if Miquel point $P$ is orthocenter or one of the points $M_{A0}$, $M_{B0}$, $M_{C0}$ of triangle $\triangle A_0 B_0 C_0$, then $\triangle A_0 B_0 C_0 \sim \triangle A_k B_k C_k$ for each $k \equiv 2 \pmod{3}$ and $k \equiv 0 \pmod{3}$.

**Corollary 4.** Let consider the queue of Miquel triangles $\{\triangle A_k B_k C_k\}$, където $k = 0, 1, 2, 3, \ldots$. The proved theorems lead to the conclusion that for the Miquel triads with fixed Miquel point is in force cyclic recurrence for the positions of this point $O \to H \to L \to O \ldots$, where $O$ is center of the circum-circle, $H$ - orthocenter, $L$ - center of the incircle or excircle, i.e. Miquel point $P$ is consistently center of the circum-circle, orthocenter and center of the incircle or excircle in the different triangles of the chain. Other cyclic recurrence is connected with points $S_{A0}$, $S_{B0}$ and $S_{C0}$. For example for point $S_{A0}$ it is $S_{A0} \to M_{A1} \to Q_{A2} \to S_{A3}....$.

Valcho Milchev
Maths teacher
"Petko Rachov Slaveikov" Secondary School
Kardzhali, Bulgaria

Address: 54 Belomorski Blvd., apt. 64, 6600 Kardzhali
Country: Bulgaria
e-mail: milchev_v@abv.bg
milchev.vi@gmail.com